\documentclass[titlepage,draft,12pt]{article} 
\usepackage{amsfonts,latexsym} 
\usepackage{pstricks}
\usepackage[a4paper]{geometry}

\date{}

\newcommand{\Ehat}{\widehat{E}}
\newcommand{\Ghat}{\widehat{G}}
\newcommand{\Ehep}{\widehat{E}_{\ep}}

\newcommand{\ep}{\varepsilon}
\newcommand{\qed}{{\penalty 10000\mbox{$\quad\Box$}}}
\newcommand{\re}{\mathbb{R}}

\newtheorem{thm}{Theorem}[section]

\newtheorem{rmk}[thm]{Remark}
\newtheorem{prop}[thm]{Proposition}

\title{Optimal decay estimates for the general solution to a class of 
semi-linear dissipative hyperbolic equations}

\author{Marina Ghisi\vspace{1ex}\\ 
{\normalsize Universit\`a degli Studi di Pisa} \\
{\normalsize Dipartimento di Matematica}\\ 
{\normalsize PISA (Italy)}\\
{\normalsize e-mail: \texttt{ghisi@dm.unipi.it}}
\and
Massimo Gobbino\vspace{1ex}\\ 
{\normalsize Universit\`a degli Studi di Pisa} \\
{\normalsize Dipartimento di Matematica}\\ 
{\normalsize PISA (Italy)}\\  
{\normalsize e-mail: \texttt{m.gobbino@dma.unipi.it}}
\and
Alain Haraux\vspace{1ex}\\ 
{\normalsize Universit\'{e} Pierre et Marie Curie} \\
{\normalsize Laboratoire Jacques-Louis Lions}\\ 
{\normalsize PARIS (France)}\\  
{\normalsize e-mail: \texttt{haraux@ann.jussieu.fr}}}

\begin{document}
\maketitle
\begin{abstract}
	We consider a class of semi-linear dissipative hyperbolic 
	equations in which the operator associated to the linear part 
	has a nontrivial kernel.
	
	Under appropriate assumptions on the nonlinear term, we prove 
	that all solutions decay to 0, as $t\to +\infty$, at least as 
	fast as a suitable negative power of $t$. Moreover, we prove that 
	this decay rate is optimal in the sense that there exists a 
	nonempty open set of initial data for which the corresponding 
	solutions decay exactly as that negative power of $t$.
	
	Our results are stated and proved in an abstract Hilbert space 
	setting, and then applied to partial differential equations.
	
\vspace{1cm}

\noindent{\bf Mathematics Subject Classification 2010 (MSC2010):}
 35B40, 35L71, 35L90.

% 35B40: Asymptotic behavior of solutions
% 35L71: Semilinear second-order hyperbolic equations
% 35L90: Abstract hyperbolic equations

\vspace{1cm} 

\end{abstract}

%%%%%%%%%%%%%%%%%%%%%
%                   %
%   Inizio lavoro   %
%                   %
%%%%%%%%%%%%%%%%%%%%%
 
\section{Introduction} 

The present work has its origin in the search for decay estimates of
solutions to some evolution equations of the general form
\begin{equation}
	u''(t)+u'(t)+Au(t)+f(u(t))=0,
	\label{eqn:f(u)}
\end{equation}
where $H$ is a real Hilbert space, $A$ is a nonnegative self-adjoint
linear operator on $H$ with dense domain, and $f$ is a
nonlinearity tangent to $0$ at the origin. \bigskip

When $f\equiv 0$, then for rather general classes of strongly positive
operators $A$ it is known that all solutions decay to $0$ (as
$t\rightarrow +\infty$) exponentially in the energy norm.  Therefore,
by perturbation theory it is reasonable to
expect that also all solutions of (\ref{eqn:f(u)}) which decay to 0
have an exponential decay rate.  The situation is different when $A$
has a non-trivial kernel.  In this case solutions tend to 0 if $f$
fulfils suitable sign conditions, but we do not expect all solutions
to have an exponential decay rate.  Let us consider for example the
hyperbolic equation
\begin{equation}
	u_{tt}+u_{t}-\Delta u+|u|^{p}u=0,
	\label{eqn0:neumann}
\end{equation}
with homogeneous Neumann boundary conditions in a bounded domain
$\Omega$.  In \cite{dec-est}, by relying on the so-called
\L ojasiewicz gradient inequality \cite{Loja1,Loja2},
it was established that, for any sufficiently small integer $p$, all
solutions of this problem tend to $0$ in the energy norm at least as
fast as $t^{-1/p}$.  Showing the optimality of this estimate means
exhibiting a ``slow solution'', namely a solution decaying exactly as
$t^{-1/p}$. \bigskip 

The existence of slow solutions for the Neumann problem was proved
in~\cite{dec-est} in the special case $p = 2$.  The main idea is that
each solution $v(t)$ to the ordinary differential equation
\begin{equation}
	v''+ v' + |v|^{p}v=0
	\label{eqn:ODE}
\end{equation}
corresponds to the spatially homogeneous solution $u(t, x) := v(t)$ of
(\ref{eqn0:neumann}), so that it is enough to exhibit a family of
solutions of (\ref{eqn:ODE}) decaying exactly as $t^{-1/2}$.  It was
later shown in \cite{Slow-fast} that actually any solution of
(\ref{eqn0:neumann}) tends to $0$ either exponentially or exactly as
$t^{-1/p}$.  This is the so-called ``slow-fast alternative''.
Moreover at this occasion the set of initial data producing
exponentially decaying solutions was shown to be closed with empty
interior.  In particular the set of ``slow" solutions corresponds to
an open set of initial data, but apart from the spatially homogeneous
solutions no explicit condition on the initial data was found in
\cite{Slow-fast}.

The proofs of these results seem to exploit in an essential way the
fact that the kernel of the linear part (in this case the set of
constant functions) is an invariant space for (\ref{eqn0:neumann}).
Without this assumption, both the alternative and the optimality of
decay rates remained open problems.

Indeed let us consider, as a model case, the hyperbolic equation
\begin{equation}
	u_{tt}+u_{t}-\Delta u-\lambda_{1}u
	+|u|^{p}u=0
	\label{eqn0:dirichlet}
\end{equation}
with homogeneous Dirichlet boundary conditions (here $\lambda_{1}$
denotes the first eigenvalue of $-\Delta$ in $H_0^1(\Omega)$).  Now
the kernel of the operator is the first eigenspace, which is not
invariant by the nonlinear term, and even the existence of a slow
solution decaying exactly as $t^{-1/p}$ was unknown until now.\bigskip

In this paper we consider a general evolution equation of type
(\ref{eqn:f(u)}), with $f$ a gradient operator satisfying some
regularity and structure conditions.  Our aim is twofold.  To begin
with, in Theorem~\ref{thm:main-above} we establish a general upper
estimate of the energy, valid for all solutions.  This estimate is
proved in a quite general context through a modified Lyapunov
functional, without any analyticity assumption on $f$.  Then in
Theorem~\ref{thm:main-below} we prove the existence of slow solutions.
This is the main result of this paper.

Our abstract theory applies to both (\ref{eqn0:neumann}) and
(\ref{eqn0:dirichlet}).  This shows in particular that the natural
upper energy estimate for solutions of these problems is in general
optimal, thereby settling an open problem raised in \cite{dec-est} and
not solved, even for the special case (\ref{eqn0:dirichlet}), by the
results of \cite{Slow-fast}.\bigskip

The problem of slow solutions has already been considered in the
parabolic setting, and in particular in the case of equation
\begin{equation} 
	u_{t}-\Delta u+|u|^{p}u=0
	\label{eqn0:neumann-parab}
\end{equation} 
with homogeneous Neumann boundary conditions in a
bounded domain $\Omega$, and in the case of equation
\begin{equation}
	u_{t}-\Delta u-\lambda_{1}u
	+|u|^{p}u=0
	\label{eqn0:dirichlet-parab}
\end{equation}
with homogeneous Dirichlet boundary conditions.  In the case of
(\ref{eqn0:neumann-parab}), an easy application of the maximum
principle shows that all solutions decay to 0 in $L^\infty (\Omega)$
at least as fast as $t^{-1/p}$ as $t\rightarrow +\infty$.  The same
property is true for (\ref{eqn0:dirichlet-parab}) but more delicate to
establish (see for example~\cite{HJK}).  With Neumann boundary conditions,
the optimality of this decay rate can be confirmed by looking at
spatially homogenous solutions as in the hyperbolic setting.  With
Dirichlet boundary conditions, a comparison with suitable
sub-solutions proves that all solutions with nonnegative initial data
are actually slow solutions (see~\cite{HJK} for the details), which
verifies the optimality of the upper estimate also in this second
case.  Moreover, in the case of Neumann boundary conditions, the
slow-fast alternative is known (see~\cite{BA}), fast solutions are
known to be ``exceptional'', and some explicit classes of slow
solutions with a sign changing initial datum were found in
\cite{BA-H}.  On the contrary, in the case of Dirichlet boundary
conditions, even the slow-fast alternative is presently an open
problem.\bigskip

All results for these parabolic problems rely on the existence of
special invariant sets, or on comparison arguments.  Both tools do not
extend easily to second order equations of the general
form~(\ref{eqn:f(u)}).  For this reason, in this paper we follow a
different path.  The main idea is to look for slow solutions in the
place where they are more likely to be, namely close to the kernel of
$A$.  Thus, under the assumption that $|f(u)|\sim|u|^{p+1}$, we look
for solutions of~(\ref{eqn:f(u)}) such that
\begin{equation}
	\langle Au(t),u(t)\rangle\leq C|u(t)|^{2p+2}
	\quad\quad
	\forall t\geq 0
	\label{est:main}
\end{equation}
for a suitable constant $C$.  Roughly speaking, under this condition
the term $Au(t)$ in~(\ref{eqn:f(u)}) can be neglected, and the dynamical 
behavior is decided by the nonlinearity only.  Thus we are in a situation
analogous to the ordinary differential equation (\ref{eqn:ODE}), for
which the existence of slow solutions can be easily established.  In
order to prove (\ref{est:main}), one is naturally led to consider the
quotient
$$Q_{p}(t):=\frac{\langle Au(t),u(t)\rangle}{|u(t)|^{2p+2}},$$
which seems to be a $p$-extension of the Dirichlet quotient (the same 
quantity with $p=0$), well known in many questions concerning 
parabolic problems (see for example the classical papers 
\cite{BT,ghidaglia} or the more recent \cite{kukavica}).

The Dirichlet quotient is nonincreasing in the case of linear
homogeneous parabolic equations.  This could naively lead to guess the
monotonicity, or at least the boundedness, of $Q_{p}(t)$ also in the
case of the second order problem (\ref{eqn:f(u)}).  Of course this is
not true as stated, but it is true for a hyperbolic version of
$Q_{p}(t)$ with a kinetic term in the numerator.  Thus we obtain the
energy $G(t)$ defined by (\ref{defn:G}), which in turn we perturb by
adding a mixing term, in such a way that the final energy $\Ghat(t)$
given by (\ref{defn:Ghat}) satisfies a reasonable differential
inequality.  This strategy is inspired by similar modified Dirichlet
quotients introduced in~\cite{ghisi:JDE2006}, and then largely
exploited in \cite{gg:k-decay,gg:w-dg} in the context of Kirchhoff
equations.  In those papers the setting is different (quasi-linear
instead of semi-linear), the goal is different
(in~\cite{ghisi:JDE2006} the main problem is the existence of global
solutions), but the strategy is the same (comparing solutions of
partial differential equations with solutions of ordinary differential
equations), thus similar tools can be applied.\bigskip

Our method produces not only some special slow solution, but an open
set in the basic energy space.  This is the first step towards proving
that slow solutions are in some sense generic, in accordance with the
general idea that the slowest decay rate is dominant, and faster
solutions are somewhat atypical.  We plan to consider this issue in a
future research.\bigskip

This paper is organized as follows.  In section~\ref{sec:statements}
we clarify the functional setting, we recall the notion of weak
solutions, and we state our main abstract results.  In
section~\ref{sec:proofs} we prove them.  In
section~\ref{sec:applications} we present some applications of our
theory to dissipative hyperbolic equations.

\setcounter{equation}{0}
\section{Functional setting and main abstract results}\label{sec:statements}

We consider the semilinear abstract second order equation
\begin{equation}
	u''(t)+u'(t)+Au(t)+\nabla F(u(t))=0
	\quad\quad
	\forall t\geq 0,
	\label{pbm:eq}
\end{equation}
with initial data
\begin{equation}
	u(0)=u_{0},
	\quad
	u'(0)=u_{1}.
	\label{pbm:data}
\end{equation}

We always assume that $H$ is a Hilbert space, and $A$ is a
self-adjoint linear operator on $H$ with dense domain $D(A)$.  We
assume that $A$ is nonnegative, namely $\langle Au,u\rangle\geq 0$ for
every $u\in D(A)$, so that for every $\alpha\geq 0$ the power
$A^{\alpha}u$ is defined provided that $u$ lies in a suitable domain
$D(A^{\alpha})$, which is itself a Hilbert space with norm
$$|u|_{D(A^{\alpha})}:=\left(|u|^{2}+|A^{\alpha}u|^{2}\right)^{1/2}.$$

We assume that $F:D(A^{1/2})\to\re$.  When we write $\nabla F(u)$, we
mean that there exists a function $\nabla F:D(A^{1/2})\to H$ such that
\begin{equation}
	\lim_{|v|_{D(A^{1/2})}\to 0}
	\frac{F(u+v)-F(u)-\langle\nabla 
	F(u),v\rangle}{|v|}=0
	\quad\quad
	\forall u\in D(A^{1/2}).
	\label{defn:gradient}
\end{equation}

The existence of $\nabla F(u)$ in the sense of (\ref{defn:gradient})
is enough to guarantee the continuity of $F$ with respect to the norm
of $D(A^{1/2})$.  Moreover, for every $u\in C^{1}([0,+\infty);H)\cap
C^{0}([0,+\infty);D(A^{1/2}))$ we have that the function $t\to
F(u(t))$ is of class $C^{1}$, and its time-derivative can be computed with
the usual chain rule $$\frac{d}{dt}\left[F(u(t))\right]= \langle\nabla
F(u(t)),u'(t)\rangle \quad\quad
\forall t\geq 0.$$

We always assume that $\nabla F:D(A^{1/2})\to H$ is locally Lipschitz
continuous, namely
\begin{equation}
	|\nabla F(u)-\nabla F(v)|\leq
	L\left(|u|_{D(A^{1/2})},|v|_{D(A^{1/2})}\right)\cdot
	|u-v|_{D(A^{1/2})}
	\label{hp:nabla-lip}
\end{equation}
for every $u$ and $v$ in $D(A^{1/2})$, for a suitable function
$L:\re^{2}\to\re$ which is bounded on bounded sets.  Under these
hypotheses, one obtains the following result concerning global
existence, regularity and derivatives of energies.

\begin{prop}\label{thmbibl:global}
	Let $H$ be a Hilbert space, let $A$ be a self-adjoint nonnegative
	operator on $H$ with dense domain $D(A)$, and let
	$F:D(A^{1/2})\to\re$.
	
	Let us assume that 
	\begin{enumerate}
		\renewcommand{\labelenumi}{(\roman{enumi})} 

		\item $F(u)\geq 0$ for every $u\in D(A^{1/2})$,
	
		\item $F$ has a gradient $\nabla F:D(A^{1/2})\to H$ in the 
		sense of (\ref{defn:gradient}),
		
		\item $\nabla F$ is locally Lipschitz continuous in the 
		sense of (\ref{hp:nabla-lip}).
	\end{enumerate}
	
	Then, for every $(u_{0},u_{1})\in D(A^{1/2})\times H$, problem
	(\ref{pbm:eq})--(\ref{pbm:data}) admits a unique global weak
	solution
	\begin{equation}
		u\in C^{0}\left([0,+\infty);D(A^{1/2})\right)\cap
		C^{1}\left([0,+\infty);H\right).
		\label{reg:w}
	\end{equation}
	
	In addition the functions 
	\begin{equation}
		E_0(t) := \frac{1}{2}\left( 
		{|u'(t)|^{2}+|A^{1/2}u(t)|^{2}}\right), 
		\quad\quad\quad
		F_0(t) := E_0(t) + F(u(t))
		\label{defn:energies}
	\end{equation}
	are of class $C^1$, and their time-derivative is given by 
	\begin{equation}
		E_0'(t) = - |u'(t)|^{2} -\langle\nabla F(u(t), u'(t)\rangle,
		\quad\quad\quad 
		F_0'(t) = -|u'(t)|^{2}.
		\label{deriv:energies}
	\end{equation}
\end{prop}\bigskip

The first main result of this paper is an upper energy estimate, valid for all 
weak solutions of (\ref{pbm:eq}).

\begin{thm}[Upper decay estimate for  weak solutions]\label{thm:main-above}
	Let us assume that 
	\begin{list}{(Hp\arabic{enumi})}{\usecounter{enumi}
		\setlength{\labelwidth}{4em}\setlength{\leftmargin}{3em}}
		
		\item $H$ is a Hilbert space, and $A$ is a self-adjoint nonnegative
		operator on $H$ with dense domain $D(A)$,
		
		\item $F:D(A^{1/2})\to[0,+\infty)$ is a nonnegative function
		with $F(0)=0$,
	
		\item $F$ has a gradient $\nabla F:D(A^{1/2})\to H$ in the 
		sense of (\ref{defn:gradient}),
		
		\item $\nabla F$ is locally Lipschitz continuous in the 
		sense of (\ref{hp:nabla-lip}),
	
		\item  there exists a constant $K>0$ such that
		\begin{equation}
			\langle\nabla F(u),u\rangle\geq K\cdot F(u)
			\quad\quad
			\forall u\in D(A^{1/2}),
			\label{hp:above-nabla}
		\end{equation}
	
		\item there exist $p>0$, and a function $R_{1}:\re\to\re$
		which is bounded on bounded sets, such that
		\begin{equation}
			|u|^{p+2}\leq R_{1}(|u|_{D(A^{1/2})})\cdot
			\left(|A^{1/2}u|^{2}+F(u)\right)
			\quad\quad
			\forall u\in D(A^{1/2}).
			\label{hp:above-p}
		\end{equation}
	\end{list}
	
	Let $(u_{0},u_{1})\in D(A^{1/2})\times H$, and let $u(t)$ be the
	unique global weak solution of problem
	(\ref{pbm:eq})--(\ref{pbm:data}) provided by
	Proposition~\ref{thmbibl:global}.
	
	Then there exist constants $M_{1}$ and $M_{2}$ such that 
	\begin{equation}
		|u'(t)|^{2}+|A^{1/2}u(t)|^{2}+F(u(t))\leq
		\frac{M_{1}}{(1+t)^{1+2/p}}
		\quad\quad
		\forall t\geq 0,
		\label{th:above-1}
	\end{equation}
	\begin{equation}
		|u(t)|\leq
		\frac{M_{2}}{(1+t)^{1/p}}
		\quad\quad
		\forall t\geq 0.
		\label{th:above-2}
	\end{equation}
\end{thm}

Our second main result is the existence of an open set of slow
solutions, namely solutions for which (\ref{th:above-2}) is optimal.

\begin{thm}[Existence of slow solutions]\label{thm:main-below}
	Let us assume that hypotheses (Hp1) through (Hp4) of 
	Theorem~\ref{thm:main-above} are satisfied. In addition, let us 
	assume that
	\begin{equation}
		\ker A\neq\{0\},
		\label{hp:A-ker}
	\end{equation}
	\begin{equation}
		\exists\nu>0\mbox{ such that $|A^{1/2}u|^{2}\geq\nu
		|u|^{2}\quad\forall u\in D(A^{1/2})\cap \ker
		(A)^{\perp}$},
		\label{hp:A-coercive}
	\end{equation} 
	and that there exist real numbers $\rho>0$, $R> 0$, $\alpha>0$ 
	such that
	\begin{equation}
		|\nabla F(u)| \leq R
		\left(|u|^{p+1}+|A^{1/2}u|^{1+\alpha}\right)
		\label{hp:below}
	\end{equation}
	for every $u\in D(A^{1/2})$ with $|u|_{D(A^{1/2})}\le \rho$.

	Then there exist a nonempty open set $\mathcal{S}\subseteq
	D(A^{1/2})\times H $ and a constant $M_{3}$ such that, for every
	$(u_{0},u_{1})\in\mathcal{S}$, the unique global solution of
	problem (\ref{pbm:eq})--(\ref{pbm:data}) provided by
	Proposition~\ref{thmbibl:global} satisfies
	\begin{equation}
		|u(t)|\geq \frac{M_{3}}{(1+t)^{1/p}} \quad\quad \forall t\geq 0.
		\label{th:below}
	\end{equation}
\end{thm}

\begin{rmk}
	\begin{em} 
		Condition (\ref{hp:A-coercive}) is known to be equivalent to
		the property that $A$ has closed range $R(A) = (\ker A)^{\perp} $.  
		
		Let $P:H\to\ker A$ denote the orthogonal
		projection on $\ker A$, and let $Q=I-P$ denote the 
		orthogonal projection on $R(A)$. From 
		(\ref{hp:A-coercive}) and (\ref{th:above-1}) it follow that
		$$|Qu(t)|^{2}\leq\frac{1}{\nu}|A^{1/2}u(t)|^{2}\leq 
		\frac{M_{1}}{\nu}\frac{1}{(1+t)^{1+2/p}}.$$
		
		Since $|u|^{2}=|Pu|^{2}+|Qu|^{2}$ for every $u\in H$, 
		comparing with (\ref{th:below}) we obtain that there exists a 
		constant $M_{4}$ such that
		$$|Pu(t)|\geq \frac{M_{4}}{(1+t)^{1/p}} 
		\quad\quad
		\forall t\geq 0.$$
		
		In other words, the range component decays faster, and the
		slow decay of $u(t)$ is due to its component with respect to
		$\ker A$.  This extends to the general abstract setting what
		previously observed in the special case studied in
		\cite{Slow-fast}.
	\end{em}
\end{rmk}

\setcounter{equation}{0}
\section{Proofs}\label{sec:proofs} 

\subsection{Proof of Proposition~\ref{thmbibl:global}}\label{sec:proof-global}

\paragraph{\textmd{\emph{Local existence}}}

We consider the Hilbert space $\mathcal{H}:=D(A^{1/2})\times H$,
endowed with the norm defined by 
$$|U|_{\mathcal{H}}^2 = |(u,v)|_{\mathcal{H}}^2 = |u|_{D(A^{1/2})}^2 
+ |v|^2,$$ 
the subspace $D(\mathcal{A}):=
D(A)\times D(A^{1/2})$, the linear operator $$\mathcal{A}(u,v):=(-v,
Au+u) \quad\quad
\forall (u,v)\in D(\mathcal{A}),$$
and the operator
$$\mathcal{F}(u,v):=( 0, u-\nabla F(u)-v)
\quad\quad
\forall (u,v)\in\mathcal{H}.$$

It is easy to check that $\mathcal{A}$ is a skew-adjoint linear
operator, hence in particular a maximal monotone linear operator on
$\mathcal{H}$ with dense domain $D(\mathcal{A})$, and
$\mathcal{F}:\mathcal{H}\to\mathcal{H}$ is a locally Lipschitz
continuous operator.  Introducing $U(t) := (u(t), u'(t))$, one can
rewrite problem (\ref{pbm:eq})--(\ref{pbm:data}) in the form $$
U'(t)+\mathcal{A}U(t)=\mathcal{F}(U(t)) \quad\quad
\forall t\geq 0,$$
with initial datum $U(0)=U_{0}:=(u_{0},u_{1})$. Thus we have reduced 
our problem to the framework of Lipschitz perturbations of maximal 
monotone operators. At this point, local existence follows from 
classical results, for which we refer to Theorem~4.3.4 and 
Proposition~4.3.9 of~\cite{C-H}. More precisely, we obtain the 
following.

\begin{itemize}
	\item \emph{(Local existence of weak solutions)} For every
	$(u_{0},u_{1})\in D(A^{1/2})\times H$, there exists $T>0$ such
	that problem (\ref{pbm:eq})--(\ref{pbm:data}) has a unique weak
	solution $$u\in C^{0}\left([0,T);D(A^{1/2})\right)\cap
	C^{1}\left([0,T);H\right).$$

	\item \emph{(Continuation)} The local solution can be continued to
	a solution defined in a maximal interval $[0,T_{*})$, with either
	$T_{*}=+\infty$, or
	$$\limsup_{t\to T_{*}^{-}}
	\left(|u'(t)|^{2}_{H}+|u(t)|^{2}_{D(A^{1/2})}\right)=+\infty.$$
	
\end{itemize}

\paragraph{\textmd{\emph{Differentiation of  energies}}}

We show that for all weak solutions the functions $E_0(t)$ and
$F_{0}(t)$ defined by (\ref{defn:energies}) are of class $C^1$, and
their time-derivative is given by (\ref{deriv:energies}) for every
$t\in[0,T)$.  Indeed for the first result we can consider the isometry
group generated on $\mathcal{H}$ by $\mathcal{A}$.  Then Lemma~11 of
\cite{LN841} (see  also \cite{Strauss} for an earlier more general
result in the same direction) gives 
$$ E_0'(t) + \langle
u(t), u'(t)\rangle =\langle\mathcal{F}(U(t)),
U(t)\rangle_{\mathcal{H}} = \langle
u(t)-\nabla F(u) -u'(t), u'(t)\rangle $$
yielding the proper result for $E_0$.  The result for $F_0$ follows
also since $\langle\nabla F(u(t)), u'(t)\rangle$ is the derivative of
the $ C^1$ function $ F(u(t))$ as a consequence of the chain rule, as
already observed.

\paragraph{\textmd{\emph{Global existence}}}

Thanks to the ``continuation'' result, all we need to show is that
$E_{0}(t)$ is bounded uniformly in time.  This follows at once from
the nonincreasing character of $F_{0}$ and our assumption that
$F(u)\geq0$.  

\subsection{A basic a priori estimate}

The next simple a priori estimate will be useful in the proof of both
main theorems.

\begin{prop}\label{prop:basic}
	Let $H$ be a Hilbert space, let $A$ be a self-adjoint nonnegative
	operator on $H$ with dense domain $D(A)$, and let
	$F:D(A^{1/2})\to\re$.
	
	Let us assume that
	\begin{enumerate}
		\renewcommand{\labelenumi}{(\roman{enumi})}
		\item  $F(u)\geq 0$ for every $u\in D(A^{1/2})$,
	
		\item $F$ has a gradient $\nabla F:D(A^{1/2})\to H$ in the 
		sense of (\ref{defn:gradient}),
	
		\item  $\langle\nabla F(u),u\rangle\geq 0$ for every $u\in 
		D(A^{1/2})$.
	
	\end{enumerate}
	
	Let $(u_{0},u_{1})\in D(A^{1/2})\times H$, and let $u(t)$ be the
	local weak solution of problem (\ref{pbm:eq})--(\ref{pbm:data}) in
	some time-interval $[0,T)$.  Then we have
	\begin{equation}
		|u'(t)|^{2}+|u(t)|^{2}+|A^{1/2}u(t)|^{2}+F(u(t))\leq
		16\left(|u_{1}|^{2}+|u_{0}|^{2}+|A^{1/2}u_{0}|^{2}+F(u_{0})\right)
		\label{th:basic-est}
	\end{equation}
	for every $t\in[0,T)$.
\end{prop}

\paragraph{\textmd{\textit{Proof}}}

Let us consider the two different energies
$$\widetilde{E}(t):=|u'(t)|^{2}+\frac{1}{2}|u(t)|^{2}+|A^{1/2}u(t)|^{2}+
	2F(u(t))+\langle u'(t),u(t)\rangle,$$
$$\Ehat(t):=|u'(t)|^{2}+|u(t)|^{2}+|A^{1/2}u(t)|^{2}+F(u(t)).$$

Due to assumption (i) and inequality
$$\left|\langle u'(t),u(t)\rangle\right|\leq\frac{3}{8}|u(t)|^{2}+
\frac{2}{3}|u'(t)|^{2},$$
it is easy to see that
\begin{equation}
	\frac{1}{8}\Ehat(t)\leq \widetilde{E}(t)\leq 2\Ehat(t)
	\quad\quad
	\forall t\in[0,T).
	\label{est:E-Ehat-basic}
\end{equation}

The function $\widetilde{E}(t)$ is of class $C^{1}$, even
in the case of weak solutions, and its time-derivative
is
$$\widetilde{E}'(t)=-|u'(t)|^{2}-|A^{1/2}u(t)|^{2}-\langle\nabla
F(u(t)),u(t)\rangle.$$

From assumption~(iii) we see that $\widetilde{E}'(t)\leq 0$, hence
$\widetilde{E}(t)\leq \widetilde{E}(0)$ for every $t\in[0,T)$.
Keeping (\ref{est:E-Ehat-basic}) into account, we have proved that
$$\Ehat(t)\leq 8\widetilde{E}(t)\leq 8\widetilde{E}(0)\leq 16\Ehat(0)
\quad\quad
\forall t\in[0,T),$$
which is exactly (\ref{th:basic-est}).\qed

\subsection{Proof of Theorem~\ref{thm:main-above}}

Let us describe the strategy of the proof before entering into 
details. We consider the energies
\begin{equation}
	E(t):=|u'(t)|^{2}+|A^{1/2}u(t)|^{2}+2F(u(t)) = 2F_0(t),
	\label{defn:E-above}
\end{equation}
\begin{equation}
	\Ehep(t):=E(t)+\ep\left[E(t)\right]^{\beta}\langle 
	u'(t),u(t)\rangle,
	\label{defn:Ehat-above}
\end{equation}
where $\ep>0$ is a parameter and
\begin{equation}
	\beta:=\frac{p}{p+2}.
	\label{defn:beta}
\end{equation}

Now we claim three facts (from now on, all positive constants
$\ep_{0}$, $\ep_{1}$, $c_{0}$, \ldots, $c_{10}$ depend on $p$,
$|u_{0}|$, $E(0)$, $K$, and on the function $R_{1}$).
\begin{itemize}
	\item  \emph{First claim}. There exist $c_{0}$ and $c_{1}$ such that
	\begin{equation}
		E(t)\leq c_{0}
		\quad\quad
		\forall t\geq 0,
		\label{est:E}
	\end{equation}
	\begin{equation}
		|u(t)|^{p+2}\leq c_{1}E(t)
		\quad\quad
		\forall t\geq 0.
		\label{est:u-E}
	\end{equation}

	\item  \emph{Second claim}. There exists $\ep_{0}>0$ such that
	\begin{equation}
		\frac{1}{2}E(t)\leq\Ehep(t)\leq 2E(t)
		\quad\quad
		\forall t\geq 0,\ \forall\ep\in(0,\ep_{0}].
		\label{est:E-Ehep}
	\end{equation}

	\item  \emph{Third claim}. There exist $\ep_{1}\in(0,\ep_{0}]$, and a 
	constant $c_{2}>0$, such that
	\begin{equation}
		\Ehep'(t)\leq -c_{2}\ep\left[\Ehep(t)\right]^{1+\beta}
		\quad\quad
		\forall t\geq 0,\ \forall\ep\in(0,\ep_{1}].
		\label{est:deriv-Ehep}
	\end{equation}
\end{itemize}

If we prove these three claims, then we easily obtain (\ref{th:above-1}) and
(\ref{th:above-2}).  Indeed let us integrate the differential
inequality (\ref{est:deriv-Ehep}) with $\ep=\ep_{1}$. We obtain the inequality
$$\widehat{E}_{\ep_{1}}(t)\leq\frac{c_{3}}{(1+t)^{1/\beta}}
\quad\quad
\forall t\geq 0.$$

Thanks to (\ref{est:E-Ehep}) and (\ref{defn:beta}), this proves
(\ref{th:above-1}).  At this point, (\ref{th:above-2}) follows from
(\ref{th:above-1}) and (\ref{est:u-E}).  So we are left to proving our
three claims.

\subparagraph{\textmd{\textit{Proof of first claim}}}

We can apply Proposition~\ref{prop:basic}. Thus
obtain estimate (\ref{est:E}) and the boundedness of
$|u(t)|_{D(A^{1/2})}$.  Then (\ref{est:u-E}) follows from
(\ref{est:E}) and assumption (\ref{hp:above-p}).

\subparagraph{\textmd{\textit{Proof of second claim}}}

From (\ref{est:u-E}) and (\ref{defn:E-above}) we have 
\begin{equation}
	\left|\langle u'(t),u(t)\rangle\right|\leq
	|u'(t)|\cdot|u(t)|\leq
	\left[E(t)\right]^{1/2}\cdot 
	c_{4}\left[E(t)\right]^{1/(p+2)},
	\label{est:scalar-prod}
\end{equation}
hence 
$$\left[E(t)\right]^{\beta}\left|\langle u'(t),u(t)\rangle\right|\leq
c_{4}\left[E(t)\right]^{p/(2p+4)}\cdot E(t).$$

Since $p>0$, with the help of (\ref{est:E}) we deduce
$$\left[E(t)\right]^{\beta}\left|\langle u'(t),u(t)\rangle\right|\leq
c_{5}E(t)\quad\quad \forall t\geq 0.$$

This implies that (\ref{est:E-Ehep}) holds true provided that 
$c_{5}\ep_{0}\leq 1/2$.

\subparagraph{\textmd{\textit{Proof of third claim}}}

The time-derivative of (\ref{defn:Ehat-above}) is
\begin{eqnarray}
	\Ehep'(t) & = & 
	-2|u'(t)|^{2} -2\ep\beta\left[E(t)\right]^{-2/(p+2)}
	\langle u'(t),u(t)\rangle|u'(t)|^{2}+
	\ep\left[E(t)\right]^{\beta}|u'(t)|^{2}
	\nonumber  \\
	\noalign{\vspace{1ex}}
	 &  & -\ep\left[E(t)\right]^{\beta}\langle u'(t),u(t)\rangle
	 -\ep\left[E(t)\right]^{\beta}\left(
	 |A^{1/2}u(t)|^{2}+\langle\nabla F(u(t)),u(t)\rangle\right)
	\nonumber  \\
	\noalign{\vspace{1ex}}
	 & =: & F_{1}+F_{2}+F_{3}+F_{4}+F_{5}.
	\label{eqn:deriv-Ehep}
\end{eqnarray}

Let us estimate separately the sum $F_{2}+F_{3}$ and the two last
terms.  First (\ref{est:scalar-prod}) implies
$$\left[E(t)\right]^{-2/(p+2)}
\left|\langle u'(t),u(t)\rangle\right|\leq
c_{4}\left[E(t)\right]^{p/(2p+4)}.$$

Then, since $p> 0$, by using (\ref{est:E}) we derive
\begin{equation}
	F_{2}+F_{3}\leq\ep|u'(t)|^{2}\left(
	c_{6}\left[E(t)\right]^{p/(2p+4)}+
	\left[E(t)\right]^{\beta}\right)\leq
	c_{7}\ep|u'(t)|^{2}.
	\label{est:I2-I3}
\end{equation}

Moreover from (\ref{est:u-E}) and (\ref{defn:beta}) we infer
$$\ep\left[E(t)\right]^{\beta}\left|\langle 
u'(t),u(t)\rangle\right|\leq
\frac{1}{2}|u'(t)|^{2}+
\frac{1}{2}\ep^{2}\left[E(t)\right]^{2\beta}|u(t)|^{2}\leq
\frac{1}{2}|u'(t)|^{2}+
c_{8}\ep^{2}\left[E(t)\right]^{2\beta+2/(p+2)}.$$

Since $2\beta+2/(p+2)=\beta+1$, this means that
\begin{equation}
	F_{4}\leq\frac{1}{2}|u'(t)|^{2}+
	c_{8}\ep^{2}\left[E(t)\right]^{\beta+1}.
	\label{est:I4}
\end{equation}

Finally, from assumption (\ref{hp:above-nabla}) and the inequality
(\ref{est:E}), we deduce
\begin{eqnarray*}
	\left[E(t)\right]^{\beta}\left( |A^{1/2}u(t)|^{2}+\langle\nabla
	F(u(t)),u(t)\rangle\right) & \geq & 
	c_{9}\left[E(t)\right]^{\beta}\left(|A^{1/2}u(t)|^{2}+
	F(u(t))\right)\\
	 & \geq & \frac{c_{9}}{2}\left[E(t)\right]^{\beta}
	 \left(E(t)-|u'(t)|^{2}\right)\\
	 & \geq & c_{10}\left[E(t)\right]^{\beta+1}
	 -c_{11}|u'(t)|^{2},
\end{eqnarray*}
hence
\begin{equation}
	F_{5}\leq -c_{10}\ep\left[E(t)\right]^{\beta+1}
	+c_{11}\ep|u'(t)|^{2}.
	\label{est:I5}
\end{equation}

Plugging (\ref{est:I2-I3}) through (\ref{est:I5}) into 
(\ref{eqn:deriv-Ehep}), we now find
$$\Ehep'(t)\leq -|u'(t)|^{2}\left(
\frac{3}{2}-c_{7}\ep-c_{11}\ep\right)+
\ep\left(c_{8}\ep-c_{10}\right)
\left[E(t)\right]^{\beta+1}.$$

If we choose $\ep_{1}\in(0,\ep_{0}]$ small enough so that
$$c_{7}\ep_{1}+c_{11}\ep_{1}\leq\frac{3}{2}
\quad\quad\mbox{and}\quad\quad
c_{8}\ep_{1}-c_{10}\leq -\frac{c_{10}}{2},$$
then (\ref{est:deriv-Ehep}) holds true with $c_{2}=c_{10}/2>0$. 
This completes the proof of Theorem~\ref{thm:main-above}\qed

\subsection{Proof of Theorem~\ref{thm:main-below}}

Let us describe the strategy of the proof before entering into
details.  Let $\nu$, $\rho$, $R$, $\alpha$ be the constants appearing
in (\ref{hp:A-coercive}) and (\ref{hp:below}).  First of all, let us
choose $\delta>0$ such that
\begin{equation}
	\delta\leq\frac{\nu}{2\nu+1}.
	\label{hp:delta}
\end{equation}

Note that this condition implies in particular that 
\begin{equation}
	\delta\leq 1
	\quad\quad
	\mbox{and}
	\quad\quad
	\delta\leq\frac{\sqrt{\nu}}{2}.
	\label{hp:delta-bis}
\end{equation}

Let $Q$ denote the orthogonal projection from $H$ to $(\ker
A)^{\perp}$.  Assuming $(u_{0},u_{1})\in D(A^{1/2})\times H$ and
$u_{0}\neq 0$, we set $$\sigma_{0}:=4
\left(|u_{1}|^{2}+|u_{0}|^{2}+|A^{1/2}u_{0}|^{2}+F(u_{0})\right)^{1/2},$$
$$\sigma_{1}:=\frac{1}{|u_{0}|^{2p+2}}\left(
\frac{1}{2}|u_{1}|^{2}+\frac{1}{2}|A^{1/2}u_{0}|^{2}+
\delta|\langle u_{1},Qu_{0}\rangle|\right)+
\frac{128R^{2}}{\delta^{2}}.$$

Let $\mathcal{S}\subseteq D(A)\times D(A^{1/2})$ be the set of 
initial data such that
\begin{equation}
	\sigma_{0}<\rho,
	\quad\quad\quad
	2\sigma_{0}^{\alpha}R<\frac{\delta}{4},
	\quad\quad\quad
	4(p+1)\sigma_{0}^{p}\sqrt{\sigma_{1}}<\frac{\delta}{32}.
	\label{hp:small}
\end{equation}

It is clear that these smallness assumptions define an open set.  This
open set is nonempty because it contains at least all pairs
$(u_{0},u_{1})$ with $u_{1}=0$ and $u_{0}\in\ker A$ with $u_{0}\neq 0$
and $|u_{0}|$ small enough.  This is the point where assumption
(\ref{hp:A-ker}) and the fact that $F(0)=0$ are essential.

Now we claim that, for every pair of initial data
$(u_{0},u_{1})\in\mathcal{S}$, the global weak solution of
(\ref{pbm:eq})--(\ref{pbm:data}) satisfies 
\begin{equation}
	u(t)\neq 0
	\quad\quad
	\forall t\geq 0, 
	\label{claim-1}
\end{equation}
and
\begin{equation}
	\frac{1}{2}\frac{|u'(t)|^{2}+|A^{1/2}u(t)|^{2}}{|u(t)|^{2p+2}}\leq 
	2\sigma_{1}
	\quad\quad
	\forall t\geq 0. 
	\label{claim-2}
\end{equation}

This is enough to prove (\ref{th:below}). Indeed, setting 
$y(t):=|u(t)|^{2}$, we observe that
\begin{equation}
	|y'(t)|=2|\langle u'(t),u(t)\rangle|\leq
	2\frac{|u'(t)|}{|u(t)|^{1+p}}\cdot|u(t)|^{2+p}\leq
	4\sqrt{\sigma_{1}}\cdot|y(t)|^{1+p/2},
	\label{est:y'}
\end{equation}
and in particular 
\begin{equation}
	y'(t)\geq -4\sqrt{\sigma_{1}}\cdot|y(t)|^{1+p/2}
	\quad\quad
	\forall t\geq 0.
	\label{diff-ineq}
\end{equation}

Since $y(0)>0$, this inequality concludes the proof.  So we are left
to prove (\ref{claim-1}) and (\ref{claim-2}). To this end, we set
\begin{equation}
	G(t):=\frac{1}{2}\frac{|u'(t)|^{2}+|A^{1/2}u(t)|^{2}}{|u(t)|^{2p+2}},
	\label{defn:G}
\end{equation}
and 
$$T:=\sup\left\{t\geq 0:\forall\tau\in[0,t],\,\,\, u(\tau)\neq 0\mbox{ and }
G(\tau)\leq 2\sigma_{1} \strut\right\}.$$

Since $u(0)\neq 0$, and $G(0)<\sigma_{1}$ (because of our definition
of $\sigma_{1}$), we have that $T>0$.  We claim that $T=+\infty$,
which is equivalent to (\ref{claim-1}) and (\ref{claim-2}).
Let us assume by contradiction that this is not the case.  Due to the
maximality of $T$, this means that either $u(T)=0$ or
$G(T)=2\sigma_{1}$.  Now we show that both choices lead to an
impossibility.

Let us set as usual $y(t):=|u(t)|^{2}$.  For every $t\in[0,T)$ we have
that $u(t)\neq 0$ and $G(t)\leq 2\sigma_{1}$.  Therefore, arguing as
in (\ref{est:y'}), we obtain that the differential inequality in
(\ref{diff-ineq}) holds true for every $t\in [0,T)$.  Since $y(0)>0$,
and $1+p/2\geq 1$, this differential inequality implies that $y(T)\neq
0$, hence $u(T)\neq 0$.

So it remains to show that $G(T)<2\sigma_{1}$. To this end, we 
introduce the perturbed energy
\begin{equation}
	\Ghat(t):=
	\frac{1}{2}\frac{|u'(t)|^{2}+|A^{1/2}u(t)|^{2}}{|u(t)|^{2p+2}}+
	\delta\frac{\langle u'(t),Qu(t)\rangle}{|u(t)|^{2p+2}}.
	\label{defn:Ghat}
\end{equation}

Due to the second condition in (\ref{hp:delta-bis}), the energy $\Ghat(t)$ 
is a small perturbation of $G(t)$ in the sense that
\begin{equation}
	\frac{1}{2}G(t)\leq\Ghat(t)\leq 2G(t)
	\quad\quad
	\forall t\in[0,T).
	\label{est:G-Ghat}
\end{equation}

The correcting term $\langle u'(t),Qu(t) \rangle$ appears frequently
when looking for boundedness or decay properties for equations whose
generator has a non-trivial kernel (see~\cite{Z} or~\cite{Range}). 

The time-derivative of $\Ghat$ is
\begin{eqnarray}
	\Ghat'(t) & = & 
	-\frac{|u'(t)|^{2}}{|u(t)|^{2p+2}}
	-\delta\frac{|A^{1/2}u(t)|^{2}}{|u(t)|^{2p+2}}
	-\frac{\langle\nabla F(u(t)),u'(t)+\delta Qu(t)\rangle}{|u(t)|^{2p+2}}
	\nonumber  \\
	 &  & 
	+\delta\frac{|Qu'(t)|^{2}-\langle u'(t),Qu(t)\rangle}{|u(t)|^{2p+2}}
	-2(p+1)\frac{\langle u'(t),u(t)\rangle}{|u(t)|^{2}}\cdot\Ghat(t)
	\nonumber  \\
	 & =: & I_{1}+\ldots+I_{5}.
	\label{eqn:G'}  
\end{eqnarray}

Let us estimate $I_{3}$, $I_{4}$, and $I_{5}$. First of all, from 
Proposition~\ref{prop:basic} we obtain that
\begin{equation}
	|u(t)|^{2}+|A^{1/2}u(t)|^{2}\leq\sigma_{0}^{2}
	\quad\quad
	\forall t\geq 0.
	\label{est:basic}
\end{equation}

Therefore, from the first smallness condition in (\ref{hp:small}) and 
assumption~(\ref{hp:below}), it follows that
$$|\nabla F(u(t))|\leq 
R\left(|u(t)|^{p+1}+|A^{1/2}u(t)|^{1+\alpha}\right)\leq
R\left(|u(t)|^{p+1}+|A^{1/2}u(t)|\cdot\sigma_{0}^{\alpha}\right),$$
hence
\begin{equation}
	\frac{|\nabla F(u(t))|}{|u(t)|^{p+1}}\leq
	R\left(1+\frac{|A^{1/2}u(t)|}{|u(t)|^{p+1}}\cdot
	\sigma_{0}^{\alpha}\right)\leq
	R\left(1+\sqrt{2G(t)}\cdot\sigma_{0}^{\alpha}\right).
	\label{est:I3-1}
\end{equation}

On the other hand, from assumption~(\ref{hp:A-coercive}) and the fact 
that $\delta\leq\sqrt{\nu}$, it follows that
$$\delta|Qu(t)|\leq\frac{\delta}{\sqrt{\nu}}|A^{1/2}u(t)|\leq
|A^{1/2}u(t)|,$$
hence
\begin{equation}
	\frac{|u'(t)|+\delta|Qu(t)|}{|u(t)|^{p+1}}\leq
	\frac{|u'(t)|+|A^{1/2}u(t)|}{|u(t)|^{p+1}}\leq
	\sqrt{2G(t)}.
	\label{est:I3-2}
\end{equation}

From (\ref{est:I3-1}) and (\ref{est:I3-2}) it follows that
$$I_{3}\leq
\frac{|\nabla F(u(t))|}{|u(t)|^{p+1}}\cdot
\frac{|u'(t)|+\delta|Qu(t)|}{|u(t)|^{p+1}}\leq
R\sqrt{2G(t)}+2R\sigma_{0}^{\alpha}G(t)\leq
\frac{4R^{2}}{\delta}+\frac{\delta}{8}G(t)+2R\sigma_{0}^{\alpha}G(t).$$

From the second smallness assumption in (\ref{hp:small}) we finally 
conclude that
\begin{equation}
	I_{3}\leq\frac{4R^{2}}{\delta}+\frac{3\delta}{8}G(t)
	\quad\quad
	\forall t\in[0,T).
	\label{est:I3-below}
\end{equation}

As for $I_{4}$, we exploit that $|Qu'(t)|\leq|u'(t)|$ 
and $|Qu(t)|\leq\nu^{-1/2}|A^{1/2}u(t)|$, hence
\begin{eqnarray*}
	|Qu'(t)|^{2}+|u'(t)|\cdot|Qu(t)| & \leq &
	|u'(t)|^{2}+\frac{1}{\sqrt{\nu}}|u'(t)|\cdot|A^{1/2}u(t)| \\
	 & \leq & \left(1+\frac{1}{2\nu}\right)|u'(t)|^{2}+
	\frac{1}{2}|A^{1/2}u(t)|^{2}.
\end{eqnarray*}

Thus from (\ref{hp:delta}) we deduce that
\begin{equation}
	I_{4}\leq\delta\frac{|Qu'(t)|^{2}+|u'(t)|\cdot|Qu(t)|}{|u(t)|^{2p+2}}
	\leq\frac{1}{2}\frac{|u'(t)|^{2}}{|u(t)|^{2p+2}}
	+\frac{\delta}{2}\frac{|A^{1/2}u(t)|^{2}}{|u(t)|^{2p+2}}
	\quad\quad
	\forall t\in[0,T).
	\label{est:I4-below}
\end{equation}

In order to estimate $I_{5}$, we exploit once again (\ref{est:basic}) 
and we obtain 
$$I_{5}\leq
2(p+1)\frac{|u'(t)|}{|u(t)|^{p+1}}\cdot|u(t)|^{p}\cdot\Ghat(t)\leq
2(p+1)\sqrt{2G(t)}\cdot\sigma_{0}^{p}\cdot\Ghat(t).$$

Since $G(t)\leq 2\sigma_{1}$ for every $t\in[0,T)$, the third 
smallness condition in (\ref{hp:small}) gives
\begin{equation}
	I_{5}\leq 4(p+1)\sqrt{\sigma_{1}}\cdot\sigma_{0}^{p}\cdot\Ghat(t)
	\leq\frac{\delta}{32}\Ghat(t)
	\quad\quad
	\forall t\in[0,T).
	\label{est:I5-below}
\end{equation}

Plugging (\ref{est:I3-below}) through (\ref{est:I5-below}) into
(\ref{eqn:G'}) we obtain 
$$\Ghat'(t)\leq -\frac{1}{2}\frac{|u'(t)|^{2}}{|u(t)|^{2p+2}}
-\frac{\delta}{2}\frac{|A^{1/2}u(t)|^{2}}{|u(t)|^{2p+2}}+
\frac{4R^{2}}{\delta}+\frac{3\delta}{8} G(t)+
\frac{\delta}{32}\Ghat(t).$$

Due to the first inequality in (\ref{hp:delta-bis}), this implies 
$$\Ghat'(t)\leq 
-\frac{\delta}{2}G(t)+
\frac{4R^{2}}{\delta}+\frac{3\delta}{8} G(t)+
\frac{\delta}{32}\Ghat(t)=
-\frac{\delta}{8}G(t)+
\frac{4R^{2}}{\delta}+
\frac{\delta}{32}\Ghat(t),$$
hence by (\ref{est:G-Ghat})
$$\Ghat'(t)\leq -\frac{\delta}{32}\Ghat(t)+\frac{4R^{2}}{\delta}
\quad\quad
\forall t\in[0,T).$$

Integrating this differential inequality we easily deduce that
\begin{equation}
	\Ghat(t)\leq\left(\Ghat(0)-\frac{128R^{2}}{\delta^{2}}\right)
	\exp\left(-\frac{\delta}{32}t\right)+
	\frac{128R^{2}}{\delta^{2}}
	\quad\quad
	\forall t\in[0,T).
	\label{est:Ghat}
\end{equation}

Since we already know that $u(T)\neq 0$, we have that $G(t)$ and 
$\Ghat(t)$ are defined and continuous at least up to $t=T$. Letting 
$t\to T^{-}$ in (\ref{est:Ghat}), and exploiting (\ref{est:G-Ghat})  
and our definition of $\sigma_{1}$, we deduce that
$$G(T)\leq 2\Ghat(T)< 2\left(
\Ghat(0)+\frac{128R^{2}}{\delta^{2}}\right)\leq 2\sigma_{1}.$$

This excludes that $G(T)=2\sigma_{1}$, thus completing the proof.\qed

\setcounter{equation}{0}
\section{Applications to partial differential equations}\label{sec:applications}

\subsection{Some equations with a local  nonlinearity of power type}

The following statement represents a bridge between the abstract
theory and partial differential equations.  Here $H$ is a space of
real valued functions, and we explicitly write $|u|_{H}$ for the norm
of the function $u\in H$ (not to be confused with the absolute value
$|u|$ of the same function).  Now the abstract assumptions on $\nabla
F$ are replaced by suitable inequalities between norms, which are
going to become Sobolev type inequalities in the concrete settings.

\begin{thm}[Semi-abstract result for local equations]\label{thm:sa-local}
	Let $\mathbb{X}$ be a set and $\mu$ be a measure in $\mathbb{X}$
	with $\mu(\mathbb{X})<+\infty$.  Let $H:=L^{2}(\mathbb{X},\mu)$,
	and let $A$ be a linear operator on $H$ with dense domain $D(A)$
	satisfying assumptions~(\ref{hp:A-ker}) and (\ref{hp:A-coercive})
	of Theorem~\ref{thm:main-below}.  Let $p>0$, and let us consider
	the second order equation
	\begin{equation}
		u''(t)+u'(t)+Au(t)+|u(t)|^{p}u(t)=0.
		\label{pbm:semi-abstract}
	\end{equation}
	
	Let us assume that
	\renewcommand{\labelenumi}{(\roman{enumi})}
	\begin{enumerate}
		\item  $D(A^{1/2})\subseteq 
		L^{2(p+1)}(\mathbb{X},\mu)$, and there exists a constant $K_{1}$ 
		such that
		\begin{equation}
			\|u\|_{L^{2(p+1)}(\mathbb{X},\mu)}\leq K_{1}|u|_{D(A^{1/2})}
			\quad\quad
			\forall u\in D(A^{1/2}),
			\label{hp:sa-1}
		\end{equation}
	
		\item  there exists a constant $K_{2}$ such that
		\begin{equation}
			\left\||u|^{p}v^{2}\right\|_{L^{1}(\mathbb{X},\mu)}\leq
			K_{2}|u|^{p}_{D(A^{1/2})}\cdot|v|_{D(A^{1/2})}\cdot|v|_{H}
			\quad\quad
			\forall (u,v)\in\left[D(A^{1/2})\right]^{2},
			\label{hp:sa-21}
		\end{equation}
		\begin{equation}
			\left\||u|^{2p}v^{2}\right\|_{L^{1}(\mathbb{X},\mu)}\leq
			K_{2}|u|^{2p}_{D(A^{1/2})}\cdot|v|_{D(A^{1/2})}^{2}
			\quad\quad
			\forall (u,v)\in\left[D(A^{1/2})\right]^{2}.
			\label{hp:sa-22}
		\end{equation}
	\end{enumerate}
	
	Then we have the following conclusions.
	\begin{enumerate}
		\renewcommand{\labelenumi}{(\arabic{enumi})} 
		
		\item \emph{(Decay for all weak solutions)} For every
		$(u_{0},u_{1})\in D(A^{1/2})\times H$, problem
		(\ref{pbm:semi-abstract}), (\ref{pbm:data}) has a unique
		global weak solution with the regularity prescribed by
		(\ref{reg:w}).  Moreover there exists a constant $M_{1}$ such that
		$$\|u(t)\|_{L^{2}(\mathbb{X},\mu)}\leq\frac{M_{1}}{(1+t)^{1/p}}
		\quad\quad
		\forall t\geq 0.$$
	
		\item \emph{(Existence of slow solutions)} There exist a
		nonempty open set $\mathcal{S}\subseteq D(A^{1/2})\times H$,
		and positive constants $M_{2}$ and $M_{3}$, with the following
		property.  For every pair of initial conditions
		$(u_{0},u_{1})\in\mathcal{S}$, the unique global solution of
		problem (\ref{pbm:semi-abstract}), (\ref{pbm:data}) satisfies
		$$\frac{M_{2}}{(1+t)^{1/p}}\leq
		\|u(t)\|_{L^{2}(\mathbb{X},\mu)}\leq\frac{M_{3}}{(1+t)^{1/p}}
		\quad\quad
		\forall t\geq 0.$$
	\end{enumerate}
\end{thm}

\paragraph{\textmd{\textit{Proof}}}

Let us set
$$F(u):=\frac{1}{p+2}\int_{\mathbb{X}}|u(x)|^{p+2}\,d\mu(x).$$

We claim that 
\begin{equation}
	[\nabla F(u)](x)=|u(x)|^{p} u(x)
	\label{defn:nabla-p}
\end{equation}
is the gradient of $F$ in the sense of (\ref{defn:gradient}), and that
all the assumptions of our abstract results
(Theorem~\ref{thm:main-above} and Theorem~\ref{thm:main-below}) are
satisfied.  All constants $c_{1}$, \ldots, $c_{8}$ in the sequel
depend only on $\mu(\mathbb{X})$, $p$, $K_{1}$, $K_{2}$, and on the
coerciveness constant $\nu$ which appears in (\ref{hp:A-coercive}).
Assumption (Hp1) is trivial, so that we can concentrate on the
remaining ones.

\subparagraph{\textmd{\textit{Verification of (Hp2)}}}

Assumption (i), and the fact that $\mu(\mathbb{X})<+\infty$, imply
the following inclusions 
\begin{equation}
	D(A^{1/2})\subseteq L^{2(p+1)}(\mathbb{X},\mu)\subseteq
	L^{p+2}(\mathbb{X},\mu)\subseteq
	L^{2}(\mathbb{X},\mu).
	\label{est:sa-sobolev}
\end{equation}

Thus $F$ is finite at least for every $u\in D(A^{1/2})$.  Moreover, it
is trivial that $F(0)=0$ and $F(u)\geq 0$ for every $u\in D(A^{1/2})$.

\subparagraph{\textmd{\textit{Verification of (Hp3)}}}

Assumption (i) implies that $\nabla F(u)$, as defined
by~(\ref{defn:nabla-p}), is in $H$ for every $u\in D(A^{1/2})$.  Now
we show that for every $u$ and $v$ in $D(A^{1/2})$ we have that
\begin{equation}
	\left|F(u+v)-F(u)-\langle\nabla F(u),v\rangle\right|_{H}\leq c_{1}
	\left(|u|^{p}_{D(A^{1/2})}+|v|^{p}_{D(A^{1/2})}\right)
	|v|_{D(A^{1/2})}\cdot |v|_{H},
	\label{est:gradient-p}
\end{equation}
which clearly implies (\ref{defn:gradient}).  To
this end, we start from the inequality
$$\left|\frac{1}{p+2}\left(|a+b|^{p+2}-|a|^{p+2}\right)-|a|^{p}ab\right|
\leq (p+1)\cdot 2^{p-1}\left(|a|^{p}+|b|^{p}\right)b^{2} \quad\quad
\forall (a,b)\in\re^{2},$$
which follows from the second order Taylor's expansion of the 
function $|\sigma|^{p+2}$. Setting $a:=u(x)$, $b:=v(x)$, and 
integrating over $\mathbb{X}$, we obtain that
\begin{equation}
	\left|F(u+v)-F(u)-\langle\nabla F(u),v\rangle\right|_{H}\leq
	c_{2}\int_{\mathbb{X}}\left(|u(x)|^{p}+|v(x)|^{p}\right)|v(x)|^{2}\,dx.
	\label{est:gradient-p-int}
\end{equation}

From (\ref{hp:sa-21}) we deduce
$$\left\|(|u|^{p}+|v|^{p})\cdot v^{2}\right\|_{L^{1}(\mathbb{X},\mu)}\leq
c_{3} \left(|u|^{p}_{D(A^{1/2})}+|v|^{p}_{D(A^{1/2})}\right)
\cdot|v|_{D(A^{1/2})}\cdot|v|_{H}.$$

Plugging this estimate into (\ref{est:gradient-p-int}), we obtain 
(\ref{est:gradient-p}).

\subparagraph{\textmd{\textit{Verification of (Hp4)}}}

We prove  for every $u$ and $v$ in $D(A^{1/2})$ the inequality
\begin{equation}
	\left|\nabla F(u)-\nabla F(v)\right|_{H}^{2}\leq
	c_{4}\left(|u|^{2p}_{D(A^{1/2})}+|v|^{2p}_{D(A^{1/2})}\right)
	|u-v|_{D(A^{1/2})}^{2},
	\label{est:gradient-lip}
\end{equation}
which implies (\ref{hp:nabla-lip}).  To this end, we
start from the inequality 
$$\left|\strut|a|^{p}a-|b|^{p}b\right|\leq
(p+1)\left(|a|^{p}+|b|^{p}\right)|a-b| \quad\quad
\forall (a,b)\in\re^{2},$$
which easily follows from the mean value theorem applied to the 
function $|\sigma|^{p}\sigma$. Setting $a:=u(x)$, $b:=v(x)$, and 
integrating over $\mathbb{X}$, we obtain that
\begin{eqnarray}
	\left|\nabla F(u)-\nabla F(v)\right|_{H}^{2} & = & 
	\int_{\mathbb{X}}\left|\strut|u(x)|^{p}u(x)-|v(x)|^{p}v(x)\right|^{2}dx
	\nonumber  \\
	 & \leq & c_{5}\int_{\mathbb{X}}
	 \left(|u(x)|^{2p}+|v(x)|^{2p}\right)|u(x)-v(x)|^{2}\,dx.
	\label{est:gradient-lip-int}
\end{eqnarray}

From (\ref{hp:sa-22}) we infer
$$\left\|(|u|^{2p}+|v|^{2p})\cdot
|u-v|^{2}\right\|_{L^{1}(\mathbb{X},\mu)}\leq c_{6}
\left(|u|^{2p}_{D(A^{1/2})}+|v|^{2p}_{D(A^{1/2})}\right)
\cdot|u-v|_{D(A^{1/2})}^{2}.$$

Plugging this estimate into (\ref{est:gradient-lip-int}), we obtain 
(\ref{est:gradient-lip}).

\subparagraph{\textmd{\textit{Verification of (Hp5)}}}

It is trivially satisfied.

\subparagraph{\textmd{\textit{Verification of (Hp6)}}}

Exploiting (\ref{est:sa-sobolev}) once again, we find 
$$|u|_{H}^{p+2}=\|u\|^{p+2}_{L^{2}(\mathbb{X},\mu)}\leq c_{7}
\|u\|^{p+2}_{L^{p+2}(\mathbb{X},\mu)}= c_{7}(p+2)F(u)$$
for every $u\in D(A^{1/2})$, which proves (\ref{hp:above-p}).
 \medskip

\subparagraph{\textmd{\textit{Verification of 
assumption~(\ref{hp:below})}}}

From (\ref{est:sa-sobolev}) we find
$$|\nabla F(u)|_{H}= |u|^{(p+1)}_{L^{2(p+1)}(\mathbb{X},\mu)}\leq
c_{8}|u|^{(p+1)}_{D(A^{1/2})}$$
for every $u\in D(A^{1/2})$, 
which proves (\ref{hp:below}) with $\alpha = p$ for any $\rho>0$.\qed
\medskip

We are finally ready to apply our theory to hyperbolic partial
differential equations.  We concentrate on the model examples
presented in the introduction.  We recall that in the Dirichlet case
even the existence of a single slow solution was an open problem.
Also in the Neumann case, where existence of slow solutions was
already known, the method of this paper gives the explicit conditions
(\ref{hp:small}) for a solution to decay slowly, conditions which were
not known before.

\begin{thm}[Neumann problem]\label{ex:neumann}
	Let $\Omega\subseteq\re^{n}$ be a bounded open set with the cone 
	property. Let $p$ be a positive exponent, with no further 
	restriction if $n\in\{1,2\}$, and $p\leq 2/(n-2)$ if $n\geq 3$.
	
	Let us consider the damped hyperbolic equation
	\begin{equation}
		u_{tt}(t,x)+u_{t}(t,x)-\Delta u(t,x)+|u(t,x)|^{p}u(t,x)=0
		\quad\quad
		\forall(t,x)\in [0,+\infty)\times\Omega,
		\label{eqn:neumann}
	\end{equation}
	with homogeneous Neumann boundary conditions
	\begin{equation}
		\frac{\partial u}{\partial n}(t,x)=0
		\quad\quad
		\forall(t,x)\in[0,+\infty)\times\partial\Omega,
		\label{bc:neumann}
	\end{equation}
	and initial data
	\begin{equation}
		u(0,x)=u_{0}(x),
		\quad
		u_{t}(0,x)=u_{1}(x)
		\quad\quad
		\forall x\in\Omega.
		\label{data:neumann}
	\end{equation}
	
	Then we have the following conclusions.
	\begin{enumerate}
		\renewcommand{\labelenumi}{(\arabic{enumi})} 
		\item \emph{(Decay for all weak solutions)} For every
		$(u_{0},u_{1})\in H^{1}(\Omega)\times L^{2}(\Omega)$, problem
		(\ref{eqn:neumann}) through (\ref{data:neumann}) has a unique
		global weak solution 
		$$u\in C^{0}\left([0,+\infty);H^{1}(\Omega)\right)\cap
		C^{1}\left([0,+\infty);L^{2}(\Omega)\right).$$
		
		Moreover there exists a constant $M_{1}$ such that
		\begin{equation}
			\|u(t)\|_{L^{2}(\Omega)}\leq\frac{M_{1}}{(1+t)^{1/p}}
			\quad\quad
			\forall t\geq 0.
			\label{th:neumann-1}
		\end{equation}
	
		\item \emph{(Existence of slow solutions)} There exist a
		nonempty open set $\mathcal{S}\subseteq H^{1}(\Omega)\times
		L^{2}(\Omega)$, and positive constants $M_{2}$ and $M_{3}$,
		with the following property.  For every pair of initial
		conditions $(u_{0},u_{1})\in\mathcal{S}$, the unique global
		weak solution of problem (\ref{eqn:neumann}) through
		(\ref{data:neumann}) satisfies 
		\begin{equation}
			\frac{M_{2}}{(1+t)^{1/p}}\leq
			\|u(t)\|_{L^{2}(\Omega)}\leq\frac{M_{3}}{(1+t)^{1/p}}
			\quad\quad
			\forall t\geq 0.
			\label{th:neumann-2}
		\end{equation}
	\end{enumerate}
\end{thm}

\begin{thm}[Dirichlet problem]\label{ex:dirichlet}
	Let $\Omega\subseteq\re^{n}$ be a bounded open set with the cone
	property, and let $\lambda_{1}$ be the first eigenvalue of
	$-\Delta$ in $\Omega$, with Dirichlet boundary conditions.  Let
	$p$ be a positive exponent, with no further restriction if
	$n\in\{1,2\}$, and $p\leq 2/(n-2)$ if $n\geq 3$.
	
	Let us consider the damped hyperbolic equation
	\begin{equation}
		u_{tt}(t,x)+u_{t}(t,x)-\Delta u(t,x)-\lambda_{1}u(t,x)
		+|u(t,x)|^{p}u(t,x)=0		
		\label{eqn:dirichlet}
	\end{equation}
	in $[0,+\infty)\times\Omega$, with homogeneous Dirichlet boundary
	conditions
	\begin{equation}
		u(t,x)=0
		\quad\quad
		\forall(t,x)\in[0,+\infty)\times\partial\Omega,
		\label{bc:dirichlet}
	\end{equation}
	and initial data (\ref{data:neumann}).  
	
	Then we have the following conclusions.
	\begin{enumerate}
		\renewcommand{\labelenumi}{(\arabic{enumi})} 
		
		\item \emph{(Decay for all weak solutions)} For every
		$(u_{0},u_{1})\in H^{1}_{0}(\Omega)\times L^{2}(\Omega)$,
		problem (\ref{eqn:dirichlet}), (\ref{bc:dirichlet}),
		(\ref{data:neumann}) has a unique global weak solution $$u\in
		C^{0}\left([0,+\infty);H^{1}_{0}(\Omega)\right)\cap
		C^{1}\left([0,+\infty);L^{2}(\Omega)\right)$$
		satisfying (\ref{th:neumann-1}).
		
		\item \emph{(Existence of slow solutions)} There exist
		a nonempty open set 
		$\mathcal{S}\subseteq 
		H^{1}_{0}(\Omega))\times L^{2}(\Omega)$, 
		and positive constants $M_{2}$ and $M_{3}$, with the following
		property.  For every pair of initial conditions
		$(u_{0},u_{1})\in\mathcal{S}$, the unique global weak solution
		of problem (\ref{eqn:dirichlet}), (\ref{bc:dirichlet}),
		(\ref{data:neumann}) satisfies (\ref{th:neumann-2}).
	\end{enumerate}
\end{thm}

\paragraph{\textmd{\textit{Proof of Theorem~\ref{ex:neumann} and 
Theorem~\ref{ex:dirichlet}}}}

We plan to apply Theorem~\ref{thm:sa-local} with $\mathbb{X}=\Omega$,
and $\mu$ equal to the Lebesgue measure on $\Omega$.  Concerning the
operator $A$, we distinguish two cases.
\begin{itemize}
	\item In the case of Theorem~\ref{ex:neumann} the operator is $Au=
	-\Delta u$ with Neumann boundary conditions, so that
	$D(A)=H^{2}(\Omega)$, $D(A^{1/2})=H^{1}(\Omega)$, and $\ker
	A\neq\{0\}$ because it consists of all (locally) constant
	functions.

	\item In the case of Theorem~\ref{ex:dirichlet} the operator is
	$Au=-\Delta u-\lambda_{1}u$ with Dirichlet boundary conditions, so
	that $D(A)=H^{2}(\Omega)\cap H^{1}_{0}(\Omega)$,
	$D(A^{1/2})=H^{1}_{0}(\Omega)$, and $\ker A\neq\{0\}$ because it
	consists of the first eigenspace of $-\Delta$.
\end{itemize}

In both cases, the norms $|u|_{H}$ and $|u|_{D(A^{1/2})}$ are
equivalent to the norms $\|u\|_{L^{2}(\Omega)}$ and
$\|u\|_{H^{1}(\Omega)}$, respectively, and the coerciveness assumption
(\ref{hp:A-coercive}) is satisfied because $\Omega$ is bounded and
eigenvalues are an increasing sequence.

Now we proceed to the verification of the assumptions of
Theorem~\ref{thm:sa-local}, which is the same in both cases.  The cone
property and the boundedness of $\Omega$ guarantee the usual Sobolev
embeddings
\begin{equation}
	H^{1}(\Omega)\subseteq L^{q}(\Omega)
	\left\{
	\begin{array}{ll}
		\forall q<+\infty & \mbox{if }n\leq 2,  \\
		\noalign{\vspace{0.5ex}}
		\forall q\leq 2^{*}=2n/(n-2) & \mbox{if }n\geq 3.
	\end{array}
	\right.
	\label{sobolev}
\end{equation}

All constants $c_{1}$, $c_{2}$, $c_{3}$ in the sequel depend only on
$p$, and on the Sobolev constants.

\subparagraph{\textmd{\textit{Verification of (\ref{hp:sa-1})}}}

It follows from (\ref{sobolev}) with $q=2(p+1)$ (note that our
assumption on $p$ is equivalent to $2(p+1)\leq 2^{*}$ if $n\geq 3$).

\subparagraph{\textmd{\textit{Verification of (\ref{hp:sa-21})}}}

Let $u$ and $v$ be in $D(A^{1/2})$. If $n\leq 2$, we apply
H\"{o}lder's inequality with three terms and exponents 4, 4, 2, and we
obtain 
$$\left\|\strut|u|^{p}\cdot v\cdot v\right\|_{L^{1}(\Omega)}\leq
\|u\|^{p}_{L^{4p}(\Omega)}\cdot\|v\|_{L^{4}(\Omega)}
\cdot\|v\|_{L^{2}(\Omega)}.$$

Thus from (\ref{sobolev}) with $q=4p$ and $q=4$ we conclude that
$$\left\|\strut|u|^{p}\cdot v\cdot v\right\|_{L^{1}(\Omega)} \leq
\|u\|^{p}_{H^{1}(\Omega)}
\cdot\|v\|_{H^{1}(\Omega)}\cdot\|v\|_{L^{2}(\Omega)},$$
which is exactly (\ref{hp:sa-21}).  If $n\geq 3$, we apply H\"{o}lder's
inequality with three terms and exponents $n$, $2^{*}$, 2, and we
obtain 
$$\left\|\strut|u|^{p}\cdot v\cdot v\right\|_{L^{1}(\Omega)} \leq
\|u\|^{p}_{L^{np}(\Omega)}
\cdot\|v\|_{L^{2^{*}}(\Omega)}\cdot\|v\|_{L^{2}(\Omega)}.$$

Thus from (\ref{sobolev}) with $q=np$ (note that $np\leq 2^{*}$) and
$q=2^{*}$ we conclude that
$$\left\|\strut|u|^{p}\cdot v\cdot
v\right\|_{L^{1}(\Omega)}\leq c_{1} \|u\|^{p}_{H^{1}(\Omega)}
\cdot\|v\|_{H^{1}(\Omega)}\cdot\|v\|_{L^{2}(\Omega)},$$
which proves (\ref{hp:sa-21}) also in the case $n\geq 3$.

\subparagraph{\textmd{\textit{Verification of (\ref{hp:sa-22})}}}

Let $u$ and $v$ be in $D(A^{1/2})$.  If $n\leq 2$, we apply
H\"{o}lder's inequality with exponents 2 and 2, and then 
(\ref{sobolev}).  We derive
$$\left\||u|^{2p}\cdot
v^{2}\right\|_{L^{1}(\Omega)}\leq
\|u\|^{2p}_{L^{4p}(\Omega)}\cdot\|v\|^{2}_{L^{4}(\Omega)}\leq
c_{2}\|u\|^{2p}_{H^{1}(\Omega)}\cdot\|v\|^{2}_{H^{1}(\Omega)},$$
which proves (\ref{hp:sa-22}) in this case.  If $n\geq 3$, we apply
H\"{o}lder's inequality with exponents $n/2$ and $n/(n-2)$, and then
(\ref{sobolev}).  Since $np\leq 2^{*}$, we find
$$\left\||u|^{2p}\cdot v^{2}\right\|_{L^{1}(\Omega)}\leq
\|u\|^{2p}_{L^{np}(\Omega)}\cdot\|v\|^{2}_{L^{2^{*}}(\Omega)}\leq
c_{3}\|u\|^{2p}_{H^{1}(\Omega)}\cdot\|v\|^{2}_{H^{1}(\Omega)},$$
which proves (\ref{hp:sa-22}) also in the case $n\geq 3$.\qed

\begin{rmk}
	\begin{em}
		For the sake of simplicity and shortness, we limited ourselves
		to the model nonlinearity $g_{p}(\sigma)=|\sigma|^{p}\sigma$.
		On the other hand, all results can be easily extended, with
		standard adjustments (such as the restriction to
		$L^{\infty}$-small initial data in low dimension), to
		equations with nonlinear terms which behave as $g_{p}(\sigma)$
		just in a neighborhood of the origin.
	\end{em}
\end{rmk}

\subsection{Some nonlocal equations involving projection operators}

The following result is suited to nonlocal partial differential
equations where a power nonlinearity is applied to some integral of
the unknown, and not to the unknown itself.

\begin{thm}[Semi-abstract result for nonlocal equations]\label{ex:sa-projection}
	Let $H$ be a Hilbert space, and let $A$ be a linear operator on
	$H$ with dense domain $D(A)$ satisfying
	assumptions~(\ref{hp:A-ker}) and (\ref{hp:A-coercive}) of
	Theorem~\ref{thm:main-below}, and such that
	\begin{equation}
		\dim(\ker A)<+\infty.
		\label{hp:dim-ker}
	\end{equation} 
	
	Let $M$ be a closed vector subspace of $H$ 
	such that 
	\begin{equation}
		M^{\perp}\cap\ker A=\{0\}, 
		\label{hp:ker-M}
	\end{equation}
	where $M^{\perp}$ denotes the space orthogonal to $M$.  Let
	$P_{M}:H\to M$ denote the orthogonal projection.
	Let $p>0$, and let us consider the second order equation
	\begin{equation}
		u''(t)+u'(t)+Au(t)+|P_{M}u(t)|^{p}P_{M}u(t)=0.
		\label{pbm:sa-projection}
	\end{equation}

	Then we have the following conclusions.
	\begin{enumerate}
		\renewcommand{\labelenumi}{(\arabic{enumi})} 
		
		\item \emph{(Decay for all weak solutions)} For every
		$(u_{0},u_{1})\in D(A^{1/2})\times H$, problem
		(\ref{pbm:sa-projection}), (\ref{pbm:data}) has a unique
		global weak solution with the regularity prescribed by
		(\ref{reg:w}).  Moreover there exists a constant $M_{1}$ such
		that 
		$$|u(t)|\leq\frac{M_{1}}{(1+t)^{1/p}} \quad\quad
		\forall t\geq 0.$$
	
		\item \emph{(Existence of slow solutions)} There exist a
		nonempty open set $\mathcal{S}\subseteq D(A^{1/2})\times H $,
		and positive constants $M_{2}$ and $M_{3}$, with the following
		property.  For every pair of initial conditions
		$(u_{0},u_{1})\in\mathcal{S}$, the unique global weak solution
		of problem (\ref{pbm:sa-projection})--(\ref{pbm:data})
		satisfies
		$$\frac{M_{2}}{(1+t)^{1/p}}\leq
		|u(t)|\leq\frac{M_{3}}{(1+t)^{1/p}} \quad\quad
		\forall t\geq 0.$$
	\end{enumerate}
\end{thm}

\paragraph{\textmd{\textit{Proof}}}

Let us set
$$F(u):=\frac{1}{p+2}|P_{M}u|^{p+2}.$$

We claim that
$$[\nabla F(u)](x)=|P_{M}u|^{p}P_{M}u$$
is the gradient of $F$ in the sense of (\ref{defn:gradient}), and that
all the assumptions of our abstract results
(Theorem~\ref{thm:main-above} and Theorem~\ref{thm:main-below}) are
satisfied.  Constants $c_{1}$, $c_{2}$, $c_{3}$ in the sequel depend
only on the operator $A$, on the subspace $M$, on $p$, and on the
coerciveness constant $\nu$ which appears in (\ref{hp:A-coercive}).

Assumptions (Hp1) and (Hp2) are trivial in this case.

Assumptions (Hp3) and (Hp4) require a completely standard verification,
based on the simple fact that the real function $|\sigma|^{p}\sigma$
is of class $C^{1}$ when $p>0$.  We omit the details for the sake of
shortness.

Assumption (Hp5)  follows from the equality
$$\langle\nabla F(u),u\rangle=|P_{M}u|^{p}\langle P_{M}u,u\rangle=
|P_{M}u|^{p}|P_{M}u|^{2}
\quad\quad
\forall u\in H.$$

We have now to verify (Hp6). This requires three steps. Let 
$P:H\to\ker A $ denote the orthogonal projection on $\ker A$. 
The first step is just observing that assumption 
(\ref{hp:A-coercive}) is equivalent to
\begin{equation}
	|A^{1/2}u|^{2}\geq\nu|u-Pu|^{2}
	\quad\quad
	\forall u\in D(A^{1/2}).
	\label{est:nu}
\end{equation}

The second step consists in proving that there exists $c_{1}>0$ 
such that
\begin{equation}
	|v|^{2}\leq c_{1}|P_{M}v|^{2}
	\quad\quad
	\forall v\in\ker A.
	\label{est:ker-coercive}
\end{equation}

To this end we set
$$c_{2}:=\min\left\{|P_{M}v|^{2}:v\in\ker A,\ |v|=1\right\},$$
and we observe that the minimum exists because of assumption
(\ref{hp:dim-ker}), and it is positive because of assumption
(\ref{hp:ker-M}).  This is enough to prove that
(\ref{est:ker-coercive}) holds true with $c_{1}=c_{2}^{-1}$.

Applying (\ref{est:ker-coercive}) with $v:=Pu$, we obtain
\begin{eqnarray*}
	|Pu|^{2} & \leq & c_{1}|P_{M}(Pu)|^{2}  \\
	 & = & c_{1}|P_{M}u-P_{M}(u-Pu)|^{2}  \\
	 & \leq & 2c_{1}|P_{M}u|^{2}+2c_{1}|P_{M}(u-Pu)|^{2}  \\
	 & \leq & 2c_{1}|P_{M}u|^{2}+2c_{1}|u-Pu|^{2}.
\end{eqnarray*}

If $u\in
D(A^{1/2})$, we can now apply (\ref{est:nu}) and conclude that
$$|u|^{2}=|Pu|^{2}+|u-Pu|^{2}\leq
2c_{1}|P_{M}u|^{2}+
(2c_{1}+1)|u-Pu|^{2}$$
$$\leq 2c_{1}|P_{M}u|^{2}+
\frac{2c_{1}+1}{\nu}|A^{1/2}u|^{2},$$
hence
$$|u|^{p+2}\leq c_{3}\left(|P_{M}u|^{p+2}+|A^{1/2}u|^{p+2}\right)
\leq c_{3}\left(p+2+|A^{1/2}u|^{p}\right)(F(u)+|A^{1/2}u|^{2}),$$
which proves (Hp6).

Finally, (\ref{hp:below}) is obviously satisfied with $R=1$,
independently of $\rho$ and $\alpha$.\qed

\medskip

We conclude with two examples of application of
Theorem~\ref{ex:sa-projection}.  In a certain sense they represents
the two extremes, namely the case where $M=H$, hence as large as
possible, and the case where $M$ is one-dimensional.  We omit the
simple proofs.

\begin{thm}\label{ex:phi-n-k}
	Let $\Omega\subseteq\re^{n}$ and $p$ be as in
	Theorem~\ref{ex:neumann}.  Let us consider the
	integro-differential damped hyperbolic equation
	$$ u_{tt}(t,x)+u_{t}(t,x)-\Delta u(t,x)+
	\left(\int_{\Omega}u^{2}(t,x)\,dx\right)^{p/2}u(t,x)=0,$$
	in $[0,+\infty)\times\Omega$, with Neumann boundary
	conditions (\ref{bc:neumann}), and initial data
	(\ref{data:neumann}).
	
	Then we have the same conclusions as those of Theorem~\ref{ex:neumann}.
\end{thm}
	
	\begin{thm}\label{ex:phi-n}
		Let $\Omega\subseteq\re^{n}$ and $p$ be as in
		Theorem~\ref{ex:neumann}. Let $\{\Omega_{i}\}_{i\in I}$ be 
		the set of all connected components of $\Omega$, and let
		$\varphi\in H^{1}(\Omega)$ be a function such that
		\begin{equation}
			\int_{\Omega_{i}}\varphi(x)\,dx\neq 0
			\quad\quad
			\forall i\in I.
			\label{hp:non-orth}
		\end{equation}

		Let $g_{p}:\re\to\re$ be defined by
		$g_{p}(\sigma):=|\sigma|^{p}\sigma$ for every $\sigma\in\re$,
		and let us consider the integro-differential damped hyperbolic
		equation
		$$u_{tt}(t,x)+u_{t}(t,x)-\Delta u(t,x)+
		g_{p}\left(\int_{\Omega}u(t,x)\varphi(x)\,dx\right)\varphi(x)=0,$$
		in $[0,+\infty)\times\Omega$, with Neumann boundary conditions
		(\ref{bc:neumann}), and initial data (\ref{data:neumann}).
		
		Then we have the same conclusions as those of Theorem~\ref{ex:neumann}.

\end{thm}

One can state similar results also for the Dirichlet problem, namely
by replacing the nonlinear term in Theorem~\ref{ex:dirichlet} with the
nonlinear terms appearing in Theorem~\ref{ex:phi-n-k} or
Theorem~\ref{ex:phi-n}.  The only difference is that in the Dirichlet
case the non-orthogonality condition~(\ref{hp:non-orth}) becomes
$$\int_{\Omega}\varphi(x)e(x)\,dx\neq 0$$
for every nonzero function $e(x)$ in the first eigenspace of the 
Dirichlet Laplacian.

\label{NumeroPagine}

\end{document}